\documentclass[a4paper,12pt]{amsart}
\usepackage[usenames,dvipsnames]{xcolor}
\usepackage{hyperref, tikz-cd}
\hypersetup{colorlinks=true,linkcolor={Brown},citecolor={Brown},urlcolor={Brown}}

\usepackage{microtype}
\usepackage{amssymb,amsmath,amsfonts,amsthm,bm,latexsym,amscd,verbatim,url,nicefrac,stmaryrd,enumerate,colonequals,dsfont,appendix,stackrel,mathrsfs,cleveref}
\crefname{exm}{Example}{Examples}
\crefname{cor}{Corollary}{Corollaries}
\crefname{prop}{Proposition}{Propositions}
\crefname{rmk}{Remark}{Remarks}
\crefname{lem}{Lemma}{Lemmata}

\input xy
\xyoption{all}
\usepackage{float}
\usepackage{times}
\usepackage{graphicx}
\usepackage{fullpage,adjustbox}
\usepackage{array}

\makeatletter
\newcommand{\oset}[3][0ex]{%
	\mathrel{\mathop{#3}\limits^{
			\vbox to#1{\kern-2\ex@
				\hbox{$\scriptstyle#2$}\vss}}}}
\makeatother

\RequirePackage{xcolor}

\numberwithin{equation}{section}
\newtheorem{thm}[equation]{Theorem}
\newtheorem{prop}[equation]{Proposition}
\newtheorem{lemma}[equation]{Lemma}

\theoremstyle{definition}
\newtheorem{dfn}[equation]{Definition}

\newtheorem{rmk}[equation]{Remark}
\newtheorem{exm}[equation]{Example}
\newtheorem{assu}[equation]{Assumption}

\usepackage{etoolbox}

\newcommand{\A}{\mathbb{A}}

\newcommand{\C}{\mathbb{C}}

\newcommand{\Q}{\mathbb{Q}}
\newcommand{\bbQ}{\mathbb{Q}}

\newcommand{\R}{\mathbb{R}}

\newcommand{\Z}{\mathbb{Z}}

\newcommand{\mcC}{\mathcal{C}}

\newcommand{\mcD}{\mathcal{D}}

\newcommand{\mcO}{\mathcal{O}}
\newcommand{\mcP}{\mathcal{P}}
\newcommand{\mcR}{\mathcal{R}}

\newcommand{\mfp}{\mathfrak{p}}

\newcommand{\mfq}{\mathfrak{q}}

\newcommand{\xto}[1]{\xrightarrow{#1}}

\DeclareMathOperator{\an}{an}

\DeclareMathOperator{\Aut}{Aut}

\newcommand{\dR}{\mathrm{dR}}
\newcommand{\Betti}{\textrm{Betti}}

\DeclareMathOperator{\et}{\acute{e}t}

\DeclareMathOperator{\GL}{GL}
\DeclareMathOperator{\Gal}{Gal}

\newcommand{\gr}{\mathrm{gr}}

\newcommand{\HK}{\mathrm{HK}}

\DeclareMathOperator{\Ind}{Ind}

\DeclareMathOperator{\map}{map}

\DeclareMathOperator{\Perf}{Perf}

\DeclareMathOperator{\rig}{rig}

\DeclareMathOperator{\DA}{{{DA}}}

\DeclareMathOperator{\RigDA}{{RigDA}}

\DeclareMathOperator{\SH}{{SH}}

\DeclareMathOperator{\HW}{H}
\newcommand{\ram}{\textrm{ram}}
 \DeclareMathOperator{\Fil}{Fil}
 \DeclareMathOperator{\crys}{crys}
 \DeclareMathOperator{\per}{per}

\usepackage{xspace}

\title{Ramified periods and field of definition}
\author{Giuseppe Ancona, Drago\c{s} Fr\u{a}\c{t}il\u{a}, Alberto Vezzani}

\begin{document}

\begin{abstract}
  Let $L/K$ be an extension of number fields that is ramified above $p$.  We give a new obstruction to the  descent to $K$ of smooth projective varieties defined over   $L$. The obstruction is a matrix of $p$-adic numbers that we call ``ramified periods'' arising from the comparison isomorphism between de Rham cohomology and crystalline cohomology. 
    As an application, we give simple examples of hyperelliptic curves over $\Q(\sqrt p)$ that are isomorphic to their Galois conjugates but such that their Jacobians do not descend to $\Q$ even up to isogeny.
\end{abstract}

\maketitle
\tableofcontents
\section{Introduction}\label{sect intro}
By a classical theorem of Shimura,  we know that there are hyperelliptic curves over $\C$ which are isomorphic to their complex conjugates but such that their Jacobians do not descend to $\R$ \cite[Thm. 3]{shimura1-fieldrat}. The equations are pretty complicated and need the presence of transcendental numbers in their coefficients.
Over number fields, similar results have been obtained in genus $2$, see  \cite{mestre-genus2}, \cite[Propositions 8.2, 8.4]{Fite}.

Our main  contribution to such questions is the following. 
\begin{thm}\label{main intro}
    Let $g$ be a natural number and $p$ be a prime number such that both are congruent to 1 modulo 4 and $p$ does not divide $g+1$.
    For $a\in\Q(\sqrt p)^{\times}$, let $\mathcal{C}_a$ be  the genus $g$ hyperelliptic curve over $\Q(\sqrt p)$ defined by the affine equation
     \begin{equation}\label{eq:intro-curveCa}
         \mathcal{C}_a : y^2=x^{2g+2}-a.
     \end{equation}
    Then, for infinitely many $a\in \Q(\sqrt p)^{\times}$,   $\mathcal{C}_a$ is isomorphic to its Galois conjugate $\mathcal{C}_{\bar a}$, but there is no abelian variety $A$ over $\Q$ such that the Jacobian $J(\mathcal{C}_a)$ is isogenous to $A\times_\Q \Q(\sqrt p)$.
\end{thm}
The elements $a$ are constructed from solutions to the negative Pell equation $x^2-py^2=-1$ and can be computed in practice.
The pleasant side of this theorem is that the equation is particularly simple and the genus can be arbitrary big.

\begin{rmk}\label{R:cp to descent field of moduli}
   A related, but quite different, question   concerns the descent to the field of moduli. This has a more geometric rather than arithmetic flavor.
    Given a finite extension $K\subset L$ and a variety $X$ over $L$, one can take the base change $\overline X:=X_{\overline K}$ of $X$ to an algebraic closure and ask whether $\overline X$ descends to $K$.
    It is possible that $X$ does not descend to $K$ but $\overline X$ does.
    For example, it is a classical fact that an elliptic curve over $\overline\Q$ descends to its field of moduli $\Q(j_E)$. In particular, if the $j$-invariant is rational, the curve descends to $\Q$.  However, it is not true that an elliptic curve defined over a number field descends to $\Q$ if its $j$-invariant is rational.
    In the case of hyperelliptic curves, Huggins \cite{huggins2007fields}, has proved that if $\Aut(X)/\langle\iota\rangle$ is not cyclic, then the hyperelliptic curve descends to its field of moduli.
    Recently, it was shown in \cite{bresciani2023field} that any plane curve of degree coprime to $6$ descends to its field of moduli.
\end{rmk}

\subsection*{Idea of the proof} We pick $a\in\Q(\sqrt p)$ for which $\mathcal{C}_a$ has good reduction modulo $p$. In particular, one can make sense of the crystalline cohomology of its reduction modulo $p$, which is, in this setting, a cohomology with coefficients in $\Q_p$. 
Berthelot constructed a comparison theorem with de Rham cohomology \cite{Berthelot}:
        \begin{equation}\label{eq:intro-comp-berth}        
        \HW_{\dR} (\mathcal{C}_a/\Q(\sqrt p) ) \otimes_{\Q(\sqrt p)} \Q_p(\sqrt p) \simeq \HW_{\crys} (\mathcal{C}_a,\Q_p)\otimes_{\Q_p} \Q_p(\sqrt p).
        \end{equation}  
Then one can ask the following question.  

\

\noindent($\star$) {Is there a $\Q(\sqrt p)$-basis of $\HW_{\dR} (\mathcal{C}_a/\Q(\sqrt p) )$ that corresponds to a $\Q_p$-basis of $\HW_{\crys} (\mathcal{C}_a,\Q_p) $ through the identification \eqref{eq:intro-comp-berth}?}

\

If the curve descends to a curve \( \mcC' \) over \( \mathbb{Q} \) with good reduction modulo \( p \), then the answer is affirmative: one can choose a basis of \( \HW_{\dR} (\mcC'/\mathbb{Q} ) \) that is mapped, via the Berthelot comparison isomorphism, to \( \HW_{\crys} (\mcC',\mathbb{Q}_p) \), which is canonically identified with \( \HW_{\crys} (\mathcal{C}_a,\mathbb{Q}_p) \).  

Thus, the question \((\star)\) gives rise to an \emph{obstruction to descent} (with good reduction at \( p \)) in the form of a double coset in \( \GL(\mathbb{Q}_p) \backslash \GL(\mathbb{Q}_p(\sqrt{p})) / \GL(\mathbb{Q}(\sqrt{p})) \).  

We refer to such an element as a \emph{ramified period}. This terminology reflects two key aspects: first, in the context of comparison isomorphisms, the coordinates of a basis of one cohomology group with respect to another are typically called \emph{periods}; second, the ramification of the number field extension at \( p \) plays a crucial role in extracting nontrivial information from \((\star)\).  

\

This elementary observation serves as the starting point for our paper. However, it does not yet provide an obstruction to descent to $\mathbb{Q}$, since the argument just outlined relies on the assumption that 
$\mcC'$ has good reduction. To address this, we consider some Galois descent datum associated to the cohomology theories we compare.  
Ultimately, the ``\emph{obstruction to descent}'' criterion that we propose 
is rather intricate (see \Cref{T:main obstruction}). 
Nevertheless, when applied
to the special case of the curve $\mathcal{C}_a$ introduced above, it leads to the following:

\

\begin{thm}\label{obstruction intro}
For the curve $\mcC_a$ defined by \eqref{eq:intro-curveCa} consider the Berthelot comparison isomorphism
  \begin{equation}\label{eq:intro-comp-gth}
   \wedge^g \HW^1_{\dR} (\mathcal{C}_a/\Q(\sqrt p) ) \otimes_{\Q(\sqrt p)} \Q_p(\sqrt p) \simeq \wedge^g \HW^1_{\crys} (\mathcal{C}_a,\Q_p)\otimes_{\Q_p} \Q_p(\sqrt p). 
  \end{equation}
 Consider a non-zero vector $v$ in the $g$-th step $\Fil^g(\wedge^g \HW^1_{\dR} (\mathcal{C}_a/\Q(\sqrt p)))$ of the Hodge filtration 
  and assume that   there is a non-zero vector $w$ in $\wedge^g \HW^1_{\crys} (\mathcal{C}_a,\Q_p)$ such that $v$ and $w$ generate the same $\Q_p(\sqrt p)$-line (via the isomorphism \eqref{eq:intro-comp-gth}). 
  Assume also that    $\Q(\sqrt p) \cdot v \cap \Q_p\cdot w=0.$ 
 Then there is no abelian variety $A$ over $\Q$ such that the Jacobian $J(\mathcal{C}_a)$ is isogenous to $A\times_\Q \Q(\sqrt p).$
 \end{thm}

\begin{rmk}
The examples we construct are such  that \( \mathcal{C}_a \times_{\mathbb{Q}(\sqrt{p})} \mathbb{Q}_p(\sqrt{p}) \) descends to \( \mathbb{Q}_p \). This aligns with the fact that the obstruction we consider is global rather than local in nature. Therefore the use of \( \mathbb{Q}(\sqrt{p}) \)-coefficients arising from de Rham cohomology is essential to our approach.
\end{rmk}

\subsection*{Organization of the paper}
In \Cref{section HK} we recall the Hyodo--Kato comparison isomorphism, which is a generalization of the  Berthelot comparison theorem mentioned above. 
  In \Cref{sect ram per} ramified periods as sketched in the introduction are defined. They provide an obstruction to descent with good reduction. To drop the hypothesis of good reduction we   construct a filtered version of ramified periods (\Cref{sect filt ram per}). We prove the generalization of \Cref{obstruction intro} in \Cref{sect obstr}. Finally, in \Cref{sect hyperell} we apply this obstruction to hyperelliptic curves and prove \Cref{main intro}.

\subsection*{Acknowledgements}  We thank Emiliano Ambrosi, Olivier De Gaay Fortman and Francesc Fit\'e for useful discussions.

This project started on the  island of Nisyros, during a summer school on periods. We thank Bruno Klingler and Yiannis Sakellaridis for the wonderful organization and the warm welcome.

This research was partially supported by the grants CYCLADES ANR--23--CE40--0011, DAG-ARTS  ANR-24-CE40-4098 of \emph{Agence National de la Recherche}, by the PNRR Grant CF 44/14.11.2022 and by
  the PRIN Grant 20222B24AY by \emph{Ministero dell'Universit\`a e della Ricerca}.

\section{The Hyodo--Kato isomorphism}\label{section HK}

Throughout the article we will work with the following convention.

\begin{assu}\label{conventions}
    We fix $L/K $  an extension of number fields, and we fix a prime ideal $\mfp$ of $L$. We let $\mfq$ [resp. $p$] be the prime ideal of $K$ [resp. of $\Q$] over which it lies. We also let $\hat{L}=L_{\mfp}$ be the $\mfp$-adic completion of $L$ and $\hat{K}$ be the $\mfq$-adic completion of $K$. They are finite extensions of $\Q_p$ with residue field $l$, respectively $k$. We let $\hat{L}_0$ be the field  $\hat L_0\colonequals W_{\mcO(\hat{K})}(l)[1/p]$ which is the  maximal unramified extension of  $\hat{K}$ in $\hat{L}$.
\end{assu}

We recall generalities on the Hyodo--Kato isomorphism and the motivic setting where it is stated, following \cite{deg-niz} and \cite{bgv}.  

\subsection{Recollections on the Hyodo--Kato cohomology and isomorphism}
There are two natural Weil cohomology theories for algebraic (even: rigid analytic) varieties over $\hat{L}$: the \emph{de Rham} cohomology and the \emph{Hyodo--Kato} cohomology. They both satisfy \'etale descent and are $\A^1$-invariant; in other words, they are \emph{motivic}. We briefly recall what that means and entails, and how they compare to each other.

\begin{dfn}
    We let $\DA({L})\colonequals\SH_{\et}(L,\Q)$  be the category of rational, \'etale motives over ${L}$ in the sense of e.g. \cite[D\'efinition 4.5.21]{ayoub-t2}, \cite[\S 2.1]{ayoub-weil}.  It is a compactly generated, symmetric monoidal category. For any variety $X/L$, we let $M(X)\in\DA(L)$ be the associated (homological) motive.
\end{dfn}

\begin{rmk}
Similarly, 
one can define the rigid analytic  motivic category   $\RigDA(\hat{L})$ and there is an analytification functor 
\[\DA(L)\to\RigDA(\hat{L}), \qquad M\mapsto M^{\an}\]
that is monoidal, preserves Tate-twists and sends a set of compact generators to a set of compact generators (see \cite[\S 2.1]{agv}, \cite[Proposition 2.31]{ayoub-weil}). \end{rmk}

\begin{rmk}
    Any field extension $L'/L$ induces a (left adjoint) monoidal functor $\DA(L)\to\DA(L')$ which sends the motive of a smooth variety $X$ to the motive of the variety $X_{L'}$ obtained by base change. We typically denote this functor by $M\mapsto M\times_L{L'}$ or simply by $M\mapsto M_{L'}$.
\end{rmk}

\begin{dfn}\label{def good reduction}
    We say that a motive $M\in \DA(L)$ \emph{has good reduction} (at $\mfp$) if its analytification $M^{\an}$ lies in the subcategory $\RigDA_{\gr}(\hat{L})$ which is generated (under colimits, shifts and twists) by analytic varieties with good reduction, i.e., admitting a smooth formal model over $\mcO_{\hat{L}}$ (see \cite[Definition 4.5]{bgv}). The resulting full subcategory is denoted by $\DA_\gr(L)$.
\end{dfn}

\begin{rmk}
    By definition, all proper varieties with good reduction at $\mfp$ lie in $\DA_{\gr}(L)$. More remarkably, all proper varieties with \emph{poly-stable }reduction at $\mfp$ (e.g. with semi-stable reduction) also lie in $\DA_{\gr}(L)$, see \cite[Proposition 3.29]{bkv}.
\end{rmk}

\begin{dfn}
    The de Rham homology functor 
    from smooth quasi-compact and quasi-separated varieties over $L$ to (the derived category of) $L$-vector spaces extends canonically to a monoidal functor \[\mcR_{\dR}\colon \DA(L)\to\mcD( L)\]
    which we call the \emph{de Rham homology} realization.  The dual groups of its homology will be denoted, as usual, by $\HW_{\dR}^*(M/L)$ or simply by $\HW_\dR^*(M)$ if we do not need to emphasize the coefficients. 
    Considering the Hodge filtration, this functor can be extended to a functor \[\DA(L)\to \Fil \mcD(L).\]
    See e.g. \cite[\S 4.15]{deg-niz}.

   If we let $\mcR_{\rig}\colon \DA(l)\to \mcD(\hat{L}_0)$ be the rigid realization functor (with coefficients in $\hat{L}_0$) then we can define the Hyodo--Kato realization, which is the  monoidal functor
    \[\mcR_{\HK}\colon \DA_{\gr}(L)\to\mcD( \hat{L}_0)\]
  defined by  $X\mapsto \mcR_{\rig}(\Psi X^{\an})$ (see \cite[\S 4.6]{bgv}). The associated cohomology groups will be denoted by $\HW_{\HK}^*(X)$.   See also \cite{hk,bei,cn,ey,deg-niz,bkv}.
    \end{dfn}

    \begin{rmk}
    \label{rmk:phiN!}
    The Hyodo--Kato realization is a motivic extension of the crystalline realization which exists for proper varieties $X$ with good reduction at $\mfp$ (in this case, $\Psi M(X)$ is the motive of the special fiber). It  can actually be extended  to a monoidal functor 
    \[\mcR_{\HK}\colon \DA_{\gr}(L)\to\mcD_{(\varphi,N)}(\hat L_0)\]
    where the category on the right denotes the derived category of $(\varphi,N)$-modules over $\hat L_0$. We will not use this extra structure in what follows.
    \end{rmk}

Also the comparison between crystalline and de Rham cohomology extends to the Hyodo--Kato case as follows (see the references above and \cite[Corollary 4.54]{bgv}).
\begin{thm}[The Hyodo--Kato isomorphism]
\label{T:Berthelot-local}
    There is an equivalence  between $\mcR_{\dR}$ and $\mcR_{\HK}\otimes_{\hat{L}_0}\hat{L}$, as functors from $\DA_{\gr}(\hat{L})$ to $\mcD(\hat{L})$. In particular, for any motive $M$ in $\DA_{\gr}(L)$, there is a functorial equivalence 
    \[
    \HW^*_{\dR}(M)\otimes_L\hat{L}\simeq \HW^*_{\HK}(M)\otimes_{\hat{L}_0}\hat{L}.
    \]
\end{thm}

\begin{rmk}\label{rmk:trivial_per}
Note that the left hand side of the equivalence is canonically equipped with its de Rham filtration, while the right hand side of the equivalence has the structure of a $(\varphi,N)$-module over $\hat{L}_0$. As such, the Hyodo--Kato isomorphism induces a ``syntomic'' (or ``ramified'') homological realization (cfr. \cite{deg-niz})
    \[
   \mcR_{\ram}\colon   \DA_{\gr}(L)\to \Fil\mcD(L)\times_{\mcD(\hat{L})}\mcD_{\varphi,N}(\hat L_0)
    \]
    or, forgetting some of the extra structures:
    \[
    \mcR_{\ram}\colon     \DA_{\gr}(L)\to \Fil\mcD(L)\times_{\mcD(\hat{L})}\mcD(\hat{L}_0)
    \qquad
    \mcR_{\ram}\colon     \DA_{\gr}(L)\to \mcD(L)\times_{\mcD(\hat{L})}\mcD(\hat{L}_0).
    \]
Note that in the case $L=K$, we obviously have $\hat{L}_0=\hat{K}$, and the realization $\mcR_{\ram}$ coincides with the de Rham realization. 
\end{rmk}

\begin{exm}\label{eg:dR&dR}
    Pick a motive $M\in \DA_{\gr}(\hat{L}_0)$. The Hyodo--Kato isomorphism gives an isomorphism of $\hat{L}_0$-vector spaces
    \[
    \HW^*_{\dR}(M)\simeq \HW^*_{\HK}(M).
    \]
    Note also that $\Psi M\simeq \Psi (M_{\hat{L}})$ as objects in $\DA(l)$. As such, we obtain natural isomorphisms of $\hat{L}_0$-vector spaces
    \[
    \HW^*_{\dR}(M)\simeq \HW^*_{\HK}(M)\simeq  \HW^*_{\HK}(M_{\hat{L}}).
    \]
    
    Suppose now one has also an object $N\in\DA_{\gr}(L)$ and an isomorphism $N_{\hat{L}}\xto{\sim}M_{\hat{L}}$, which induces the following commutative square:
    \[
\xymatrix{
\HW^*_{\HK}(M_{\hat{L}})\otimes_{\hat{L}_0}\hat{L}\ar[r]_-{\sim}^-{\hat{L}_0}\ar[d]_{\sim}^{\hat{L}_0}&  \HW^*_{\dR}(M/\hat L_0)\otimes_{\hat L_0}\hat L\ar[d]^-{\sim}\\
\HW^*_{\HK}(N_{\hat{L}})\otimes_{\hat{L}_0}\hat{L}\ar[r]^-{\sim}&  \HW^*_{\dR}(N/L)\otimes_L{\hat{L}})
    }
    \]
  Using what shown above, the top isomorphism (and not just the left one) descends as an isomorphism between $\hat{L}_0$-vector spaces. In particular, the Hyodo--Kato isomorphism for $N_{\hat{L}}$ can be interpreted as an isomorphism between  de Rham cohomologies of the two models  $N$ and $M$ of $N_{\hat{L}}$, each one carrying a model over $L$ resp. $\hat{L}_0$.
\end{exm}

\subsection{The Galois equivariant version of the Hyodo--Kato isomorphism}\label{sect Gal HK}
We give now the  Galois-equivariant version of the Hyodo--Kato isomorphism from \Cref{T:Berthelot-local}.
This version is needed to deal with varieties (or, more generally, motives) having \emph{potentially} good reduction and it will be used in \Cref{sect obstr}.

\begin{dfn}
    Let $\DA_{\hat{L}\gr}(\hat{K})$ [resp. $\DA_{\hat{L}\gr}({K})$] be the full subcategory of $\DA(\hat{K})$  [resp. of $\DA({K})$] whose objects $M$ are such that $M_{{\hat L}}$ lies in $\DA_{\gr}(\hat{L})$. 
\end{dfn}
\begin{rmk}
    The category  $\DA_{\hat{L}\gr}(\hat{K})$ contains $\DA_{\gr}(\hat{K})$. By \'etale descent, assuming for simplicity that $\hat{L}/\hat{K}$ is Galois, the former can be described as the category $\DA_{\gr}(\hat{L})^{\Gal(\hat{L}/\hat{K})}$ of motives with good reduction with a Galois-descent datum over $\hat{K}$ (note that $\DA_{\gr}(\hat{L})$ is Galois-stable in $\DA(\hat{L})$).

\end{rmk}

\begin{dfn}
    Assume for simplicity that $\hat{L}/\hat{K}$ is Galois. We may define the Hyodo--Kato realization on $ \DA_{\hat{L}\gr}(\hat{K})$ as follows:
       \[
    \mcR_{\HK^+}\colon \DA_{\hat{L}\gr}(\hat{K})\simeq\DA_{\gr}(\hat{L})^{\Gal(\hat{L}/\hat{K})}\to \mcD(\hat{L}_0)^{\Gal(\hat{L}/\hat{K})}
    \]
   landing in the
    category of $\Gal(\hat{L}/\hat{K})$-semilinear representations.  
     Note that the base change along $\hat{L}_0\to\hat{L}$ defines a functor  \[\mcD(\hat{L}_0)^{\Gal(\hat{L}/\hat{K})}\to \mcD(\hat{L})^{\Gal(\hat{L}/\hat{K})}\simeq \mcD(\hat{K}).\]
\end{dfn}

\begin{rmk}
    As in \Cref{rmk:phiN!}, we can actually enrich the functor $\mcR_{\HK^+}$ as a functor to the category  $\mcD_{(\varphi,N)}(\hat{L}_0)^{\Gal(\hat{L}/\hat{K})}$.
We won't use this enrichment in what follows.
\end{rmk}
We have then the following form of the Hyodo--Kato isomorphism, cfr. \cite{deg-niz}.
\begin{prop}[Galois-equivariant Hyodo--Kato isomorphism]\label{P:equiv HK} 

    Assume  that $\hat{L}/\hat{K}$ is Galois. There is a  natural transformation between the realization functors
    \[
    \mcR_{\dR}\colon \DA_{\hat{L}\gr}({K})\to \Fil\mcD(K)\to \mcD(\hat{K})
    \]
    and
    \[
    \mcR_{\HK^+}\colon \DA_{\hat{L}\gr}({K})\to 
    \mcD(\hat{L}_0)^{\Gal(\hat{L}/\hat{K})}\to \mcD(\hat{K})
    \]
    giving rise to a functor
    \[
    \mcR_{per}\colon  \DA_{\hat{L}\gr}({K})\to \Fil\mcD(K)\times_{\mcD(\hat{K})} \mcD(\hat{L}_0)^{\Gal(\hat{L}/\hat{K})}.
    \]
    In particular, for any $M\in \DA_{\hat{L}\gr}({K})$ there is a $\Gal(\hat{L}/\hat{K})$-semilinear equivariant isomorphism
    \[
    \HW^*_{\HK}(M_L)\otimes_{\hat{L}_0}\hat{L}\simeq \HW^*_{\dR}(M)\otimes_K\hat{L}.
    \]
    
\end{prop}
\begin{proof}
    The functor $\RigDA_{\gr}(\hat{L})$ has descent with respect to unramified extensions of finite valued fields above $\hat{K}$ (by \cite[Corollary 4.14]{bgv}) and its \'etale sheafification coincides with $\RigDA(-)$ (this is because $\RigDA(C)=\RigDA_{\gr}(C)$ if $C$ is algebraically closed), see \cite[Theorems 3.3.3(2), 3.7.21]{agv}. Also the functor $\mcD((-)_0)$  has descent with respect to unramified extensions, and its \'etale sheafification is, by construction, the functor  $\mcD({K}^{nr})^{\Gal(-)}$ where we let $\mcD({K}^{nr})$ be $\Ind\varinjlim_F\Perf(F)$ as $F$ runs through finite unramified extensions of $K$ in $\bar{K}$,  and where $\Gal(K)$ acts via its quotient to $\Gal(k)$.    This implies that the \'etale sheafification of the Hyodo--Kato realization induces a functor
    \[
    \RigDA(\hat{K})\to \mcD({K}^{nr})^{\Gal(K)}
    \]
    which, by construction, restricts (via analytification) to the functor $\mcR_{\HK^+}$ on the category $\DA_{\hat{L}\gr}({K})$. 

    After such sheafification, the Hyodo--Kato isomorphism extends canonically to a Galois-equivariant equivalence between $\mcR_{\dR}$ and $\mcR_{\HK^+}$. When restricted to the category $\DA_{\hat{L}\gr}({K})$, it gives  the content of the statement.
\end{proof}

\section{Ramified periods}\label{sect ram per}
In this section, we define ramified periods based on the Hyodo--Kato comparison isomorphism. We show  that they form some obstruction to descending a motive  to a smaller field in the case of good reduction.

We keep the notation from \Cref{conventions} and we consider a motive $M$ over $L$ with good reduction at $\mathfrak{p}$ (in the sense of  \Cref{def good reduction}). As an example, we may consider  (the motive of) a smooth projective variety having good reduction at $\mathfrak{p}$.

\begin{dfn}\label{rem ram}\label{def ram}
Let $M$ be a motive in $\DA_{\gr}(L)$ and $n\in\Z$. 
    Consider the identification 
         \[  \HW^n_{\dR} (M/L)\otimes_L \hat{L} \simeq \HW^n_{\HK} (M,\hat{L}_0)\otimes_{\hat{L}_0} \hat{L} \]
from \Cref{T:Berthelot-local}. Let $d$ be the dimension of the vector spaces in this identification. If one chooses an $L$-basis of $\HW^n_{\dR} (M/L)$ and an $\hat{L}_0$-basis of $\HW^n_{\HK} (M,\hat{L}_0)$, one can write the coordinates of one basis with respect to the other and get a matrix in $\GL_d(\hat{L})$. 
Accounting for the choice of bases, we obtain a well-defined double coset
 \[
    \mathcal{P}_{\ram}(M,n)\in\GL_d(L) \backslash \GL_d(\hat{L})/ \GL_d(\hat{L}_0). 
 \]
 We call it the \emph{$n$-th ramified period of $M$}.
\end{dfn}

\begin{rmk}\label{R:periods as double coset}
Based on the Betti--de Rham comparison theorem for varieties over $\Q$
\[
    \HW^i_{\dR}(X/\Q)\otimes_{\Q}\C\simeq \HW^i_\Betti(X(\C),\Q)\otimes\C
\] 
classical periods are defined as the matrix (with entries in $\C$) of the isomorphism taken with respect to some choice of rational bases on both sides. 
This matrix is only well defined up to multiplication on the left and on the right by invertible rational matrices.
In this setting, one typically considers the   \emph{ring of periods} $\mcP_\C(X)$ as the $\Q$-subalgebra of $\C$ generated by all such entries.
This is a well defined object independent on the choice of bases.

In our $p$-adic setting, the analogous ring would not be interesting: 
since  the compositum $\hat{L}_0\cdot L$ coincides with the whole $\hat{L}$, 
the ring generated by the entries of $\mcP_{\ram}(M,n)$ would always be the whole of $\hat{L}$. 
Even in the classical context, considering the period as a double coset in $\GL_d(\Q)\backslash \GL_d(\C)/\GL_d(\Q)$ is a finer invariant that has never been considered in the literature (as far as we are aware). In our setting this point of view is needed in order to have a non trivial information. Albeit elementary, this is a crucial observation of this work.
\end{rmk}

Let us first see that there is a lot of space to have non-trivial periods:
\begin{lemma}\label{L: a lot of periods} 
If $\hat{L}/\hat{K}$ is ramified, then for all $d\ge 1$
     the double quotient 
     \[
     \GL_d(L) \backslash \GL_d(\hat{L})/ \GL_d( \hat{L}_0)
     \]
     has cardinality of the continuum.
\end{lemma}
   \begin{proof}
       If $r $ and $s$ are the dimensions of $\hat{L}$ and $\hat{L}_0$ as $\Q_p$-vector spaces, then $\GL_d(\hat{L})$ and $\GL_d(\hat{L}_0)$ are $p$-adic Lie groups of dimension $d^2r$ and $d^2s$. In particular, the quotient  $\GL_d(\hat{L})/ \GL_d(\hat{L}_0)$  is a $p$-adic analytic variety of dimension $d^2(r-s),$ which is positive because $r>s$ by the hypothesis. 
       One concludes since $\GL_d(L) $ is countable.
   \end{proof} 
 
\begin{rmk}
   The space $\GL_d(L) \backslash \GL_d(\hat{L})/ \GL_d(\hat{L}_0)$ does not have a nice  topology. As the proof above shows, it has the shape of something like $\Q_p/\Q$. The only information we will use in this paper is whether an element is trivial or not in this double quotient. We do not know yet if finer information can be read off of it.
\end{rmk}

\begin{rmk}
    The target category of the realization 
    \[
    \mcR_{per}\colon \DA_{\gr}(L)\to \mcD(L)\times_{\mcD(\hat{L})}\mcD(\hat{L}_0)
    \] 
    of \Cref{rmk:trivial_per}  can be written as the category of triples $(V,W,\alpha)$ with $V\in\mcD(L)$, $W\in\mcD(
    \hat{L}_0)$ and $\alpha\colon V_{\widehat{L}}\simeq W_{\widehat{L}}$, with   mapping spaces  computed by the fiber of
    \[
    \map(V,V')\oplus\map(W,W')\stackrel{\alpha^*-\alpha'_*}{\longrightarrow}\map(V_{\hat{L}},W'_{\hat{L}'}).
    \] 
    By trivializing the complexes $V$ and $W$, any object $(V,W,\alpha)$ is given by a direct sum of shifts of elements $\alpha$ in $\GL(\widehat{K})$ and $\alpha\simeq\beta$ iff $\alpha=A\beta B$ for some $A\in\GL(K)$, $B\in\GL(\widehat{K}_0)$.  These matrices correspond to the ramified periods introduced before. 
\end{rmk}

    The following shows that ramified periods are an obstruction to descent with good reduction:
 \begin{prop}\label{ostruzione buona riduzione}
   Let $M$ be a motive in   $\DA_{\gr}(L)$. If it has a model $M_0$ in $\DA_{\gr}(K)$, then the ramified periods $\mathcal{P}_{\ram}(M,n)$ are all trivial. 
   
\end{prop}
\begin{proof}
We note that the essential image of the natural functor $P\mapsto P_L=(P_L,P_{\hat{L}_0})$ from $\mcD(K)\to \mcD(L)\times_{\mcD(\hat{L})}\mcD(\hat{L}_0)$ 
    lies in the subcategory of matrices $\alpha$ which are trivial in $\GL(L)\backslash \GL(\hat{L})/\GL(\hat{L}_0)$. By functoriality, the object $\mcR_{\per}(M_{0L})$ coincides with $\mcR(M_0)_L$ where now $\mcR(M_0)$ lies in $\mcD(K)\times_{\mcD(\hat{K})}\mcD(\hat{K})=\mcD(K)$. 
 \end{proof}

\begin{rmk}
 One would like to say that if some ramified period of a motive $M$ with good reduction is non-trivial, then $M$ can not have a model $M_0$ over $K$ possibly without good reduction. This is not straightforward and a version of this criterion will be proved in \Cref{sect obstr} based on the invariants defined in \Cref{sect filt ram per}.
\end{rmk}
 
\section{Filtered ramified   periods}\label{sect filt ram per}
In this section we refine ramified periods by taking into account the Hodge filtration. This will improve \Cref{ostruzione buona riduzione}, at least in some special (but crucial) cases.

We keep notation from \Cref{conventions}. Unless otherwise stated, $M$ denotes a motive over the number field $L$ with good reduction at $\mathfrak{p}$ (in the sense of \Cref{def good reduction}). 

 \begin{dfn}
     Let $V$ be a finite dimensional  $\hat L_0$-vector space. Let $\Fil^*$ be a decreasing filtration by $\hat L$-subvector spaces in $V\otimes_{\hat L_0} \hat L.$ We say that $\Fil^*$ \emph{descends to $V$} if, for all $i$, we have
     \[ (\Fil^i \cap V)\otimes_{\hat L_0} \hat L = \Fil^i.\]
 \end{dfn}
Our definition is motivated by the following observation:
 \begin{prop}\label{P:ostruzione fil}
 For a motive $M\in \DA_\gr(\hat L)$ and a cohomological degree $n$, consider the comparison isomorphism
         \[  \HW^n_{\dR} (M/\hat L)  \simeq \HW^n_{\HK} (M,\hat L_0)\otimes_{\hat L_0} \hat L \]
 from \Cref{T:Berthelot-local}.
     Assume that there exists a motive $N_0\in\DA_\gr(\hat L_0)$, such that $M\simeq N_0 \times_{\hat L_0} \hat L$. Then the Hodge filtration descends to the Hyodo-Kato cohomology (and coincides with the Hodge filtration on $\HW^n_{\dR}(N_0/\hat L_0)$).
     In particular, the induced filtration on $\HW^n_\HK(M,\hat{L}_0)$ is independent of $N_0$.
 \end{prop}
 \begin{proof}
      By \Cref{T:Berthelot-local} applied to both $M/\hat{L}$ and $N_0/\hat{L}_0$, we have a commutative diagram of isomorphisms 
     \[
        \begin{tikzcd}
        \HW_{\dR}^n(M/\hat{L})\ar[r,"a"]\ar[d,"b"] & \HW_\HK^n(M,\hat{L}_0)\otimes_{\hat{L}_0}\hat{L}\ar[d,"c"]\\
            \HW_{\dR}^n(N_0/\hat{L}_0)\otimes_{\hat{L}_0}\hat{L}\ar[r,"d"] &\HW_\HK^n(N_0,\hat{L}_0)\otimes_{\hat{L}_0}\hat{L}
        \end{tikzcd}
     \] 
     (see also \Cref{eg:dR&dR}).
  On the one hand, the  map $b$ is compatible with the Hodge filtrations.  On the other hand, the  maps $c$ and $d$ identify the $\hat{L}_0$-structures. Hence, the Hodge filtration on $\HW_{\dR}^n(N_0/\hat{L}_0)$ defines a filtration on $\HW_\HK^n(M,\hat{L}_0)$ which induces the Hodge filtration on $\HW_{\dR}^n(M/\hat{L})$.
    \end{proof} 

  Note that, contrary to \Cref{ostruzione buona riduzione}, the criterion given above is local (at $\mfp$) and gives an obstruction to the existence of a model (with good reduction) over $\hat{L}_0$. 
  
\begin{dfn}\label{def parab hk}\label{def parab dr}     Assume that there exists a  motive with good reduction $N_0\in\DA_{\gr}(\hat{L}_0)$ such that $M_{\hat{L}}\simeq (N_0)_{\hat{L}}$. 
  We can consider the filtration on $\HW^n_{\HK} (M,\hat{L}_0)$ induced by \Cref{P:ostruzione fil}, and define $P_{\HK}$ to be the parabolic group of automorphisms of $\HW^n_{\HK} (M,\hat{L}_0)$ respecting this filtration. Similarly, $P_{\dR}$ is the parabolic group of automorphisms of $\HW^n_{\dR} (M/L)$ respecting the Hodge filtration.
\end{dfn}
\begin{rmk}
    In general, the Hyodo--Kato cohomology (of a motive over $\hat{L}$) is \emph{not} canonically equipped with a Hodge filtration. In the above definition, it is crucial to have an \emph{algebraic} model over $\hat{L}_0$ with good reduction.
\end{rmk}

 \begin{rmk}\label{rem ram fil}
 Under the assumption of \Cref{def parab dr}, 
    the identification \Cref{T:Berthelot-local}
         \[  \HW^n_{\dR} (M/L)\otimes_L \hat{L} \simeq \HW^n_{\HK} (M,\hat{L}_0)\otimes_{\hat{L}_0} \hat{L}. \]
 induces an identification of algebraic groups
         \begin{equation}\label{eq:P_dR=P_HK}
            P_{\dR}\times_L \hat{L} \simeq P_{\HK}\times_{\hat{L}_0} \hat{L}.
         \end{equation}
\end{rmk}

\begin{dfn}\label{def ram fil}
    We keep the assumptions of \Cref{def parab dr}, and we fix $n\in\Z$.   
If one chooses an $L$-basis of $\HW^n_{\dR} (M/L)$ adapted to the Hodge filtration and a $\hat{L}_0$-basis of $\HW^n_{\HK} (M,\hat{L}_0)$ adapted to the  filtration induced by \Cref{P:ostruzione fil}, one can write the coordinates
of one basis with respect to the other and get a matrix in $ P_{\dR}  (\hat{L})= P_{\HK} (\hat{L})$ (by \eqref{eq:P_dR=P_HK}). 
In particular, the   identification \eqref{T:Berthelot-local} gives a well-defined element 
    \[  \mathcal{P}_{\Fil}(M,n)\in  P_{\dR}  (L) \backslash P_{\dR}  (\hat{L}) /  P_{\HK} (\hat{L}_0) \]
    which does not depend on the choice of the bases.
\end{dfn}

As in \Cref{ostruzione buona riduzione} and \Cref{P:ostruzione fil}, we can easily deduce the following obstruction.
 \begin{prop}\label{ostruzione 2}
  Under the assumptions of \Cref{def parab dr}, if $M\in\DA_\gr(L)$ has a model $M_0$ in $\DA_{\gr}(K)$, then the ramified periods $\mathcal{P}_{\Fil}(M,n)$ are all trivial. \qed
\end{prop}
This is not yet the obstruction we want, as $M_0$ is still assumed to have good reduction. We will overcome this in \Cref{sect obstr}.

\section{An obstruction to descent}\label{sect obstr}
 In this section we improve \Cref{ostruzione buona riduzione} and \Cref{ostruzione 2}, at least in some special (but crucial) cases. The obstruction to descent we prove is based on filtered  ramified periods (\Cref{def ram fil}).

We keep notation from \Cref{conventions} and fix a motive $M$ over the number field $L$ with good reduction modulo $\mathfrak{p}$ (in the sense of \Cref{def good reduction}), such as the motive of a smooth variety with good reduction at $\mfp$.
\begin{dfn}\label{def minimal}
We say that the motive $M$ is \emph{minimal in cohomological degree $n$}

if the smallest nonzero subspace of the Hodge filtration on $\HW_\dR^n(M)$ is of dimension one.
\end{dfn}
\begin{exm}
    If $A$ is an abelian variety of dimension $g$, the motive $H^g(A)$ is of dimension $(2g)!/(g!)^2$ and is minimal (in cohomological degree $g$).
\end{exm}
    \begin{dfn}\label{def per minimal}
Assume that the motive $M$ is minimal (in cohomological degree $n$) and that there exists a motive $N_0$ defined over $\hat{L}_0$ with good reduction such that 
\[M \times_L \hat{L}\simeq N_0\times_{\hat L_0}\hat{L}.\]
We  define the \emph{minimal period} of $M$ (in cohomological degree $n$) as the element
    \[  \mathcal{P}_{\min}(M,n) \in  L^\times \backslash  \hat{L}^\times /   \hat{L}_0^\times  \]
induced by $\mathcal{P}_{\Fil}(M)$ by restricting to the smallest nonzero subspace of the Hodge filtration.
\end{dfn}

\begin{thm}\label{T:main obstruction}Suppose that the extension $L/K$ is Galois and totally ramified at $\mfq$. 
Let $M\in\DA_{\gr}(L)$ be minimal in cohomological degree $n$, and assume the following:

\begin{enumerate}
     \item\label{i:exists N_0} There exists a motive $N_0\in \DA_{\gr}(\hat{K})$ and an isomorphism $M \times_L \hat{L}\simeq N_0\times_{\hat{K}} \hat{L}$.
    \item There exists a motive $M_0\in \DA(K)$ and an isomorphism  $M \simeq M_0 \times_K L$.
\end{enumerate}
Then  the minimal period $\mathcal{P}_{\min }(M,n)$ of $M$ is  trivial.
\end{thm}

\begin{rmk}
    We point out that the second  condition does \emph{not} imply the first: indeed $M_0$ may not have good reduction, whereas $N_0$ does. Note also that the existence of an \emph{analytic} model of $M^{\an}$ in $\RigDA_{\gr}(\hat{K})$ is automatic from the equivalences $M^{\an}\in\RigDA_{\gr}(\hat{L})\simeq\RigDA_{\gr}(\hat{K})$. 
\end{rmk}
\begin{proof}

Let $G$ be the Galois group of the extension $L/K$. It coincides with the decomposition group $\Gal(\hat{L}/\hat{K})$. Consider the comparison  isomorphism 
        \begin{equation}\label{eq:comp-proof-main-minimal}   
            \HW^n_{\dR} (M/L)\otimes_L \hat{L} \simeq \HW^n_{\HK} (M,{\hat{K}})\otimes_{\hat{K}} \hat{L}
        \end{equation}     
        from   \Cref{T:Berthelot-local}.
The existence of $M_0$ implies that both sides of \eqref{eq:comp-proof-main-minimal} are equipped with a $G$-semilinear action compatible with the isomorphism; see \Cref{P:equiv HK}. 
Both sides are equipped with a filtration (see \cref{P:ostruzione fil}) and the isomorphism is filtered.

Since the Galois action respects the Hodge filtration and the Hyodo--Kato $\hat K$-structure, the isomorphism  \eqref{eq:comp-proof-main-minimal} induces $G$-equivariant isomorphisms for all $i\in \Z$
\begin{equation}\label{eq:comp-filtered-proof-mainthm}
    \Fil^i(\HW^n_{\dR} (M/L))\otimes_L \hat{L} \simeq \Fil^i (\HW^n_{\HK} (M,\hat{K}))\otimes_{\hat{K}} \hat{L}. 
\end{equation}

In particular, since $M$ is minimal in cohomological degree $n$, by considering the smallest nonzero piece of the Hodge filtration, we obtain a pair of lines $\ell_{\dR} \subset \HW^n_{\dR} (M,L)$ and   $\ell_{\HK} \subset \HW^n_{\HK} (M,\hat{K})$ together with a $G$-equivariant isomorphism
          \[  \ell_{\dR} \otimes_L \hat{L} = \ell_{\HK} \otimes_{\hat{K}} \hat{L}. \]

On the one hand, by Galois descent for vector spaces, there is a $K$-vector space $T$ of dimension one endowed with the trivial action of $G$ such that $\ell_{\dR} =T \otimes_K L$ as $G$-representations.
On the other hand, there is a character $\chi$ of $G$ such that the isotypical component
$(\ell_{\HK} \otimes_{\hat{K}} \hat{L})^{G,\chi}$ is $\ell_{\HK}$. Hence, a vector in $(T \otimes_K L)^{G,\chi} \subset \ell_{\dR}$ is sent to $\ell_{\HK}$, which proves that the minimal period $\mathcal{P}_{\min }(M,n)$ is trivial. 
\end{proof}

 \section{Hyperelliptic curves}\label{sect hyperell}
In this section we give examples of hyperelliptic curves over $K=\bbQ(\sqrt p)$ that are isomorphic to their conjugate and such that their Jacobian does not descend to $\bbQ$. The obstruction to descent will be provided by \Cref{T:main obstruction}.

The Galois conjugate of $a\in\bbQ(\sqrt p)$ will be denoted by $\bar{a}$.
For any nonzero $t$ in some field $k$ we write  $\mathcal{C}_t$ for the hyperelliptic curve over $k$ of genus $g$ defined by the affine equation
    \[ \mathcal{C}_t : y^2=x^{2g+2}-t.\]

\begin{thm}\label{thm hyperell}
    Let $g$ be a natural number and $p$ be a prime number such that both are congruent to 1 modulo 4 and $p$ does not divide $g+1$.
 Let $v\in \Q(\sqrt p)$ be such that its norm $N_{\Q(\sqrt p)/\Q}(v)$ equals $-1$. Let $b\in \Q(\sqrt p)$ be such that $v^{g+1}=b/\bar{b}$. Define $a\in \Q(\sqrt p)$ as an element of $p$-adic valuation zero and of the form $a=b^2\cdot p^n$ for some integer $n$. Then the hyperelliptic curve $\mathcal{C}_a$  is isomorphic to its Galois conjugate  $ \mathcal{C}_{\bar a}$, but there is no abelian variety $A$ over $\Q$ such that the Jacobian $J(\mathcal{C}_a)$ is isogenous to $A\times_\Q \Q(\sqrt p).$
\end{thm}
\begin{rmk}
The equation $N_{\Q(\sqrt p)/\Q}(v) = -1$ is a Pell equation with infinitely many solutions, as $p$ is congruent to 1 modulo 4, see \cite[p. 55]{Mordell}. Furthermore, since $g$ is odd, $v^{g+1}$ has norm $1$, so by Hilbert 90,   there is an element $b\in\Q(\sqrt p)$ such that $v^{g+1}=b/\bar{b}$. 
  From the construction, we get the equality $v^{2g+2}=a/\bar{a}$, therefore, there are infinitely many possible values for $a$ to which the statement of the theorem applies. This justifies the ``infinitely many'' assertion of \Cref{main intro}.
\end{rmk}
\begin{proof}
    From the construction, we have the relation $v^{2g+2}=a/\bar{a}$. It implies that the morphism \[(x,y) \mapsto (vx,v^{g+1}y)\] induces 
    an isomorphism between $\mathcal{C}_a $ and $\mathcal{C}_{\bar a}$.
    
    To show that an abelian variety as in the statement cannot exist, we use the obstruction coming from \Cref{T:main obstruction} based on the notion of minimal period (\Cref{def per minimal}).
    First, in \Cref{P:computation}, we show that $\mathcal{C}_a \times \Q_p(\sqrt p)$ has a model over $\Q_p$ with good reduction.
    Moreover, we also compute the minimal period of the motive $M:=\wedge^g \HW^1 (\mathcal{C}_a)\simeq  \HW^g(J(\mcC_a))$. 
    Then we show that the period is non-trivial in \Cref{P:non triv} (as $g(g+1)/2$ is odd according to the hypothesis on $g$).
    Finally, we conclude that the Jacobian $J(\mcC_a)$ is not isogeneous to an abelian variety defined over $\Q$ invoking \Cref{T:main obstruction}.
\end{proof}
 
\begin{prop}\label{P:computation}
    Let $a,p,g$ be as in \Cref{thm hyperell}. 
    There exists an integer $c\in\Z$ such that $a$ and $c$ are congruent modulo $p$.

    The curve $\mathcal{C}_c$ is defined over $\Q_p$ (and even over $\Q$), it  has good reduction and satisfies
\[\mathcal{C}_c \times_{\Q_p} \Q_p(\sqrt p)\cong\mathcal{C}_a \times_{\Q(\sqrt p)} \Q_p(\sqrt p).\]
     Moreover, there is an $\alpha\in\Q_p(\sqrt p)$ such that $\alpha^{2g+2}=a/c$.
For such an $\alpha$, the minimal period (\Cref{def per minimal}) of the motive $M:=\wedge^g \HW^1 (\mathcal{C}_a)$ is given by
    \[\mathcal{P}_{\min }(M)=\alpha^{g(g+1)/2}.\]
\end{prop}
\begin{proof} 
As the residue field of $\Q_p(\sqrt p)$ is $\mathbb{F}_p$ the existence of $c$ follows. Since $a$ is nonzero modulo $p$, so is $c$ and the curve $\mathcal{C}_c$ has good reduction.

Since $a/c$ equals $1$ modulo $p$, the equation $\alpha^{2g+2}=a/c$ can be solved by Hensel's Lemma, since we assumed that $p$ and $2g+2$ are coprime.
The map 
\begin{align*}
\phi \colon \mcC_c\to\mcC_a\quad
\text{ given by } \quad (x,y) \mapsto (\alpha^{-1} x,\alpha^{-(g+1)}y)
\end{align*}
induces an isomorphism between  $\mathcal{C}_c $ and $\mathcal{C}_{a}$ over $\Q_p(\sqrt p)$. Notice that the former is defined over $\Q_p$ (and even over $\Q$) and it has good reduction as $c$ is an invertible integer modulo $p$. 
In particular, the Hyodo-Kato\footnote{same as crystalline in this case} cohomology of $\mathcal{C}_a $ is identified with the de Rham  cohomology of $\mathcal{C}_{c}$ and we can compute the period $\mathcal{P}_{\min }(M)$ through this identification.
More precisely, from \Cref{eg:dR&dR} we have a commutative square of isomorphisms
\[
\begin{tikzcd}
    \HW^1_\dR(\mcC_a/\Q(\sqrt p))\otimes\Q_p(\sqrt p) \ar[d,"\phi^*_\dR"]\ar[r,"\mcP_a"] & \HW^1_\HK(\mcC_a,\Q_p)\otimes\Q_p(\sqrt p)\ar[d,"\phi^*_\HK"]\\
    \HW^1_\dR(\mcC_c/\Q_p)\otimes\Q_p(\sqrt p) \ar[r,"\mcP_c"] & \HW^1_\HK(\mcC_c,\Q_p)\otimes\Q_p(\sqrt p)
\end{tikzcd}
\]
Since $\mcP_c$ and $\phi^*_\HK$ are defined over $\Q_p$, the minimal period $\mathcal{P}_{\min}(M)$ is the determinant of the  isomorphism induced by $\phi^*_{\dR}:$
\begin{equation}\label{eq:Fil1dR-dRisom}
    \Fil^1 \HW^1_{\dR}(\mcC_a/\Q(\sqrt p))\otimes_{\Q(\sqrt p)}\Q_p(\sqrt p) \simeq \Fil^1 \HW^1_{\dR}(\mcC_c/\Q_p)\otimes_{\Q_p}\Q_p(\sqrt p),
\end{equation}
viewed as an element in $\Q_p^\times\backslash \Q_p(\sqrt p)^\times / \Q(\sqrt p)^\times$.

A basis of $\Fil^1 \HW^1_{\dR}$ is given by the classes of the differential forms $x^i dx/y$ for $i$ between $0$ and $g-1$. 
With respect to this choice of bases, the map \eqref{eq:Fil1dR-dRisom} is diagonal with coefficients $\alpha^{g-i}$, again for $i$ between $0$ and $g-1$. Taking the determinant gives the desired formula.
\end{proof}

\begin{prop}\label{P:non triv}
    For $\alpha$ as in \Cref{P:computation} and $d$ an odd integer, the class of $\alpha^{d}$ in 
    $\Q_p^\times \backslash  \Q_p(\sqrt p)^\times / \Q(\sqrt p)^\times$ is non-trivial.
\end{prop}
\begin{proof} 
By contradiction, suppose one can write 
\[\alpha^d=\beta (x+\sqrt p y)\]
with $\beta \in \Q_p$ and $x,y \in \Q$. 
Raising to the power $2g+2$ gives the relation
\[a^d/(c^d(x+\sqrt p y)^{2g+2})=\beta^{2g+2}. \]
Since the left hand side is in $\Q(\sqrt p)$ and the right hand side in $ \Q_p$ and $\Q(\sqrt p)\cap\Q_p = \Q$, we deduce $a^d/(c^d(x+\sqrt p y)^{2g+2})\in \Q$. In particular, as  $c \in \Q$, we can write
\[a^d=rw^{2g+2}\]
for some $r\in \Q$ and $w\in\Q(\sqrt p)$. Hence, we obtain the equality
\[(a/\bar{a})^d=(w/\bar w)^{2g+2}.\]
From the hypothesis of \Cref{thm hyperell}, we have the relation $v^{2g+2}=a/\bar{a}$, which implies
\begin{equation*}\label{eq:v^d2g+2=w/wbar}
(v^d)^{2g+2}=(w/\bar w)^{2g+2}.
\end{equation*}
Since $v$ and $w$ belong to $\Q(\sqrt p)$, they are real numbers, so the above equality reduces to 
\[v^d=\pm(w/\bar w).\]

However, by construction, $v$ has norm \(-1\), and since $d$ is odd, the same holds for $v^d$. Meanwhile, any number of the form \( \pm w / \bar{w} \) has norm \( 1 \), leading to the desired contradiction.
\end{proof}

\bibliographystyle{alpha}
\bibliography{biblio}   
\end{document}